\documentclass{amsart}
\usepackage{amssymb,amscd}
\advance\hoffset by -1in
\advance\voffset by -1in

\usepackage[cp1251]{inputenc}
\usepackage[russian]{babel}
\newtheorem{theorem}{Теорема}

\newtheorem{corollary}{Следствие}

\title{О про-$p$-группах с одним определяющим соотношением}

\author{А.\,Ф.\,Красников}
\address{Омский государственный университет, г.\,Омск}
\email{phomsk@mail.ru}

\begin{document}

\maketitle

\section*{Введение}
\noindent Известная теорема о свободе Магнуса \cite{Mg} показывает, что если $F$ --- свободная группа с базой $x_1,\ldots,x_n$, $H = \textup{гр }(x_1,\ldots,x_{n-1})$, $r\in F$ --- циклически несократимое слово в образующих $x_1,\ldots,x_n$, содержащее $x_n$, $R$ --- нормальная подгруппа, порожденная в группе $F$ элементом $r$, то $H\cap R = 1$.\\
В настоящей статье показано, что условие $r$ --- {\it циклически несократимое слово в образующих $x_1,\ldots,x_n$, содержащее $x_n$}, эквивалентно условию\\
$D_n(r)\not\equiv 0 \mod{{\bf Z}[F]\cdot(R-1)}$,
где $D_n$ --- производная Фокса кольца ${\bf Z}[F]$, однозначно определяемая условиями
$D_n(x_n)=1,~D_n(x_j)=0,~\mbox{при}~n\neq j$.\\
Справедливы следующая теорема и вытекающее из этой теоремы необходимое условие выполнения теоремы о свободе для про-$p$-групп с одним определяющим соотношением:
\begin{theorem}\label{tm1_gr}
Пусть $F(n)$ --- свободная про-$p$-группа с системой образующих $\{x_1,\ldots ,x_n\}$, $K\subseteq \{1,\ldots ,n\}$, $F(K)$
 --- свободная про-$p$-группа с системой образующих $\{x_j\mid j\in K\}$; $v\in F(n)$,
$N$--- замкнутая нормальная подгруппа в $F(n)$, ${\bf Z}_p[[F(n)]]$ --- пополненное групповое кольцо группы $F(n)$; $\{\partial_1,\ldots ,\partial_n\}$ --- производные Фокса кольца ${\bf Z}_p[[F(n)]]$,
$\pi:{\bf Z}_p[[F(n)]]\rightarrow {\bf Z}_p[[F(n)/N]]$ --- гомоморфизм, индуцированный
естественным гомоморфизмом $F(n)\rightarrow F(n)/N$. Если
\begin{eqnarray}\label{end_gr}
\pi(\partial_k(v))=0,~k\in \{1,\ldots ,n\}\setminus K,
\end{eqnarray}
то $v$ принадлежит замыканию $F(K)(F(K)\cap N)^{F(n)}[N,N]$ в $F(n)$.
\end{theorem}

\begin{corollary}\label{crl_pro_p_gr}
Пусть $F(n)$ --- свободная про-$p$-группа с системой образующих $\{x_1,\ldots ,x_n\}$,
$H$ --- свободная про-$p$-группа с системой образующих $\{x_1,\ldots ,x_{n-1}\}$, $r$ --- элемент группы $F(n)$, $R$ --- замкнутая нормальная подгруппа в $F(n)$, порожденная элементом $r$, $\pi:{\bf Z}_p[[F(n)]]\rightarrow {\bf Z}_p[[F(n)/R]]$ --- гомоморфизм, индуцированный
естественным гомоморфизмом $F(n)\rightarrow F(n)/R$; $\{\partial_1,\ldots ,\partial_n\}$ --- производные Фокса алгебры ${\bf Z}_p[[F(n)]]$.
Если $H\cap R=1$, то $\pi(\partial_n(r))\neq 0$.
\end{corollary}

\noindent В \cite{Gdh} показано, что в $F(2)$ элемент $r=x_1^p[x_2,x_1^p]$ не является $p$-й степенью,
$x_1\notin R$, $x_1^p\in R$ и аналог теоремы о свободе Магнуса
для $F(2)/R$ не выполняется.\\
Нетрудно видеть, что $\pi(\partial_2(r))= \pi(-x_2^{-1}x_1^{-p}x_2x_1^p+x_1^p)= 0$.\\
Т.е.  элемент $x_1^p[x_2,x_1^p]$
не удовлетворяет необходимому условию выполнения теоремы о свободе для про-$p$-групп с одним определяющим соотношением.

\section{Некоторые свойства производных Фокса}
\noindent Пусть $F$ --- свободная группа
с базой $\{g_j \mid j\in J\}$, $N$--- нормальная подгруппа в $F$.
Обозначим через ${\bf Z}[F]$ целочисленное групповое кольцо группы
$F$. Дифференцированием кольца ${\bf Z}[F]$ называется отображение
$\partial:{\bf Z}[F]\to {\bf Z}[F]$ удовлетворяющее условиям
\begin{gather}
\partial(u+v)=\partial(u)+\partial(v),\notag\\
\partial(uv)=\partial(u)v+ \varepsilon (u)\partial(v)\notag
\end{gather}
для любых $u,~v\in {\bf Z}[F]$, где $\varepsilon$ --- гомоморфизм
тривиализации $F\to 1$, продолженный по линейности на ${\bf Z}[F]$.\\
Обозначим через $D_k~(k\in J)$ производные Фокса кольца
${\bf Z}[F]$ --- дифференцирования, однозначно определяемые
условиями:
\begin{gather}
D_j(g_j)=1,~D_k(g_j)=0,~\mbox{при}~k\neq j.\notag
\end{gather}
Для $u\in {\bf Z}[F],~f\in F,~n\in N$ имеют место формулы:
\begin{gather}
D_k(f^{-1})=-D_k(f)f^{-1},~D_k(f^{-1}nf)\equiv D_k(n)f\mod{{\bf Z}[F]\cdot (N-1)};\notag\\
u-\varepsilon (u)=\sum_{j\in J} (g_j-1)D_j(u).\label{1tm01_gr}
\end{gather}
Пусть $G$ --- группа; $A$, $B$ --- подмножества множества элементов группы $G$; $C$, $D$ --- подмножества множества элементов кольца ${\bf Z}[G]$. Обозначим через $\textup{гр}(A)$ подгруппу, порожденную $A$ в $G$, через $A^G$ --- нормальную подгруппу,
порожденную $A$ в $G$, через $\gamma_k G$ --- $k$-й член нижнего центрального ряда группы $G$. Если $x,~y$ --- элементы $G$, то положим $[x,y]=x^{-1}y^{-1}xy$, $x^y=y^{-1}xy$.
Через $AB$ обозначим множество произведений вида $ab$, где $a,~b$ пробегают соответственно элементы $A$, $B$, через $[A,B]$ --- подгруппу группы $G$, порожденную всеми $[a,b]$, $a\in A$, $b\in B$. Через $CD$ обозначим множество сумм произведений вида $cd$, где $c,~d$ пробегают соответственно элементы $C$, $D$.

\begin{theorem}\label{tm3_gr}
Пусть $F$ --- свободная группа с базой $\{g_j \mid j\in J\}$, $K\subseteq J$, $F_K$
 --- подгруппа в $F$, порожденная $\{g_j\mid j\in K\}$; $v\in F$,
$N$--- нормальная подгруппа в $F$, $1\neq d$ --- целое, $\mathfrak{N}={\bf Z}[F]\cdot(N-1)+d\cdot{\bf Z}[F]$; $D_k~(k\in J)$ --- производные Фокса кольца ${\bf Z}[F]$. Тогда и только тогда
\begin{eqnarray}\label{1tm02_gr}
D_k(v)\equiv ~0\mod{\mathfrak{N}},~k\in J\setminus K,
\end{eqnarray}
когда $v\in F_K(F_K\cap N)^F[N,N]N^d$.
\end{theorem}
\begin{proof} Непосредственно проверяется, что сравнения (\ref{1tm02_gr}) следуют из $v\in F_K(F_K\cap N)^F[N,N]N^d$.
Необходимо доказать обратное.\\
Формулы (\ref{1tm01_gr}), (\ref{1tm02_gr}) показывают, что
$({v-1} - \sum_{j\in K}(g_j-1)D_j(v))\in \mathfrak{N}$. Следовательно,
элементы $1$ и $v$ принадлежат одному правому классу смежности группы
$F$ по подгруппе $F_KN$, т.е. найдется $\hat{v}\in F_K$ такой, что $v\hat{v}^{-1}\in N$. Элемент $v\hat{v}^{-1}$ обозначим через $w$.
Предположим, что $w\notin (F_K\cap N)^F[N,N]N^d$ и приведем это предположение к противоречию.
Пусть $X=\{g_j | j\in J\}$; $u\to \bar u$ --- функция, выбирающая правые шрайеровы представители $F$ по $N$,
$S$ --- множество выбранных представителей.
Так как $N = \text{гр }(sx{\overline{sx}}^{-1}|s\in S,~x\in X)$ (доказательство
см., например, в \cite{KrMr}), то найдутся $s_ix_i{\overline{s_ix_i}}^{-1},~s_i\in S,~x_i\in
X$ и ненулевые целые $k_i$ такие, что
\begin{eqnarray*}
w\equiv (s_1x_1{\overline{s_1x_1}}^{-1})^{k_1}\ldots (s_lx_l{\overline{s_lx_l}}^{-1})^{k_l}\mod{(F_K\cap N)^F[N,N]N^d}
\end{eqnarray*}
и элемент $w$ нельзя представить по модулю $(F_K\cap N)^F[N,N]N^d$ в виде произведения степеней меньшего чем $l$ числа элементов вида
$sx{\overline{sx}}^{-1},~s\in S,~x\in X$. Через $w_i$ будем обозначать элементы $s_ix_i{\overline{s_ix_i}}^{-1}$,
$i=1,\ldots,l$. Отметим, что
\begin{eqnarray}\label{1tm03_gr}
\sum_{i=1}^l k_i D_q(w_i)\equiv 0 ~\mod{\mathfrak{N}},~q\in J\setminus K.
\end{eqnarray}
Предположим, что в $\{x_1,\ldots,x_l\}$ есть не принадлежащий $F_K$ элемент. Пусть это будет $x_1$.
Приведем это предположение к противоречию.
Выберем $q\in J\setminus K$ такое, что $x_1=g_q$.
Если $x_i=g_q$, то $D_q(w_i)= D_q(s_i){s_i}^{-1}w_i+
{\overline{s_ix_i}}^{-1}-D_q(\overline{s_ix_i}){\overline{s_ix_i}}^{-1}$. Следовательно,
\begin{eqnarray}\label{1tm04_gr}
D_q(w_i)\equiv D_q(s_i){s_i}^{-1}+
{\overline{s_ix_i}}^{-1}-D_q(\overline{s_ix_i}){\overline{s_ix_i}}^{-1}~\mod{\mathfrak{N}}.
\end{eqnarray}
Если $x_i\neq g_q$, то $D_q(w_i)= D_q(s_i){s_i}^{-1}w_i-
D_q(\overline{s_ix_i}){\overline{s_ix_i}}^{-1}$. Следовательно,
\begin{eqnarray}\label{1tm05_gr}
D_q(w_i)\equiv
D_q(s_i){s_i}^{-1}-D_q(\overline{s_ix_i}){\overline{s_ix_i}}^{-1}~\mod{\mathfrak{N}}.
\end{eqnarray}
Пусть $t\in S$ и $t=u{g_q}^\varepsilon u_1$, где $\varepsilon = \pm 1$ и слово $u_1$ не
содержит в своей записи буквы $g_q$. Тогда
\begin{eqnarray*}
D_q(t){t}^{-1}=
D_q(u){u}^{-1}+D_q({g_q}^\varepsilon)({u{g_q}^\varepsilon})^{-1},
\end{eqnarray*}
т.е. $D_q(t){t}^{-1}$ --- сумма элементов вида $D_q({g_q}^\varepsilon)({u{g_q}^\varepsilon})^{-1}$, $u{g_q}^\varepsilon\in S$.
Покажем, что
\begin{eqnarray*}
{\overline{s_1x_1}}^{-1}\not\equiv \pm
D_q({g_q}^\varepsilon)({u{g_q}^\varepsilon})^{-1}~\mod{N}.
\end{eqnarray*}
Предположим противное. Если $\varepsilon = -1$, то ${\overline{s_1x_1}}^{-1}\equiv {u}^{-1}~\mod{N}$.
Но тогда $s_1=ux_1^{-1}$ и $s_1x_1{\overline{s_1x_1}}^{-1}=1$.
Пришли к противоречию.\\
При $\varepsilon = 1$ будем иметь ${\overline{s_1x_1}}^{-1}\equiv (u{g_q})^{-1}~\mod{N}$.
Но тогда $s_1=u$ и
$\overline{s_1x_1}=s_1x_1$. Пришли к
противоречию.\\
Теперь можно утверждать, что из (\ref{1tm03_gr}), (\ref{1tm04_gr}),
(\ref{1tm05_gr}) следует существование $i$ такого, что $i\neq 1,~x_i=g_q$ и
${\overline{s_ix_i}}^{-1}={\overline{s_1x_1}}^{-1}$. Но
тогда $s_i=s_1$ и потому $w_i=w_1$. Снова пришли к противоречию.\\
Полученные противоречия показывают, что $\{x_1,\ldots,x_l\}\subset F_K$.\\
Выберем минимальное $n$ с таким свойством: найдутся $f_i\in F_K$, $0\neq m_i\in {\bf Z}$ и не кончающиеся символом из $F_K\cap X$ элементы $v_i$, $\hat{v}_i\in S$, $v_if_i\hat{v}_i^{-1}\in N$ $(i=1,\ldots,n)$ такие, что
\begin{eqnarray}\label{1tm10_gr}
w\equiv (v_1f_1\hat{v}_1^{-1})^{m_1}\ldots (v_nf_n\hat{v}_n^{-1})^{m_n}\mod{(F_K\cap N)^F[N,N]N^d}.
\end{eqnarray}
(Числом с таким свойством будет, например, $l$).
Без потери общности рассуждений мы можем и будем считать, что $v_1$ - элемент максимальной длины среди
элементов $v_1, \hat{v}_1,\ldots,v_n, \hat{v}_n$ и ему нет равных в множестве элементов $\hat{v}_1,\ldots,v_n, \hat{v}_n$.
Действительно, если $\hat{v}_1$ - элемент максимальной длины, то заменим $v_1f_1\hat{v}_1^{-1}$ на $\hat{v}_1f_1^{-1}v_1^{-1}$;
если $v_1=\hat{v}_1$, то $v_1f_1\hat{v}_1^{-1}\in (F_K\cap N)^F$ --- в противоречии с минимальностью $n$; если
$v_1=v_i$, $1\neq i$, то заменим $v_if_i\hat{v}_i^{-1}$ на $v_1f_1\hat{v}_1^{-1}\hat{v}_1f_1^{-1}f_i\hat{v}_i^{-1}$; если
$v_1=\hat{v}_i$, $1\neq i$, то заменим $v_if_i\hat{v}_i^{-1}$ на $(v_1f_1\hat{v}_1^{-1}\hat{v}_1f_1^{-1}f_i^{-1}v_i^{-1})^{-1}$.\\
Пусть $q\in J\setminus K$ такое, что $v_1$ кончается одним из символов $g_q$, $g_q^{-1}$. Из (\ref{1tm10_gr}) получаем
\begin{eqnarray}\label{1tm11_gr}
D_q(w)\equiv \sum_{i=1}^n m_i(D_q(v_i)v_i^{-1} - D_q(\hat{v}_i)\hat{v}_i^{-1})\mod{\mathfrak{N}}.
\end{eqnarray}
Если $u\in S$, то нетрудно видеть, что $D_q(u)u^{-1}$ будет суммой элементов вида $\pm t^{-1}$, $t\in S$, $t$ --- начальный отрезок слова $u$. Поэтому из (\ref{1tm11_gr}) следует, что если $v_1$ кончается символом $g_q$, то
\begin{eqnarray*}
D_q(w)\equiv m_1v_1^{-1} + \mu\mod{\mathfrak{N}},
\end{eqnarray*}
где $\mu$ --- сумма элементов вида $\pm t^{-1}$, $t\in S$, $t\neq v_1$,\\
если $v_1=\tilde{v}_1g_q^{-1}$, то
\begin{eqnarray*}
D_q(w)\equiv -m_1\tilde{v}_1^{-1} + \mu\mod{\mathfrak{N}},
\end{eqnarray*}
где $\mu$ --- сумма элементов вида $\pm t^{-1}$, $t\in S$, $t\neq \tilde{v}_1$.\\
Тогда $D_q(w)\not\equiv 0 ~\mod{\mathfrak{N}}$ --- в противоречии с (\ref{1tm03_gr}).
\end{proof}

\begin{corollary}
Пусть $F$ --- свободная группа с базой $x_1,\ldots,x_n$,
$\{D_1,\ldots ,D_n\}$ --- производные Фокса кольца ${\bf Z}[F]$,
$1\neq r$--- циклически несократимое слово в образующих $x_1,\ldots,x_n$, $R$ --- нормальная подгруппа, порожденная в группе $F$ элементом $r$, $H = \textup{гр }(x_1,\ldots,x_{n-1})$.\\
Если $H\cap R=1$, то $D_n(r)\not\equiv 0 \mod{{\bf Z}[F]\cdot(R-1)}.$
\end{corollary}
\begin{proof}
Допустим, что $H\cap R=1$, но $D_n(r)\equiv 0 \mod{{\bf Z}[F]\cdot(R-1)}.$
Теорема \ref{tm3_gr} показывает, что тогда $r\in [R,R]$.\\
Следовательно, $r\in\bigcap\limits_{n \in {\bf N}}
R^{(n)}=1$ --- противоречие.
\end{proof}
\noindent Если $D_n(r)\not\equiv 0 \mod{{\bf Z}[F]\cdot(R-1)}$, то $r\text{ содержит }x_n$ и тогда, по
теореме о свободе Магнуса, $H\cap R=1$.
Следовательно, условие $D_n(r)\not\equiv 0 \mod{{\bf Z}[F]\cdot(R-1)}$ эквивалентно условию $H\cap R=1$.

\begin{corollary}\label{clr_gr}
Пусть $F$ --- свободная группа с базой $\{g_j \mid j\in J\}$, $K\subseteq J$, $H$
 --- подгруппа в $F$, порожденная $\{g_j\mid j\in K\}$,
$R$--- нормальная подгруппа в $F$, $M$ --- множество всех тех нормальных подгрупп группы $F$, которые принадлежат $R$, $D_k~(k\in J)$ --- производные Фокса кольца ${\bf Z}[F]$.\\
Если (и только если) $H\cap R=1$, то для любых $N \in M,~v\in N$ справедлива формула
$D_k(v)\equiv 0\mod{{\bf Z}[F]\cdot(N-1)},~k\in J\setminus K\Rightarrow v\in [N,N].$
\end{corollary}
\begin{proof}
Теорема \ref{tm3_gr} показывает, что если $H\cap R=1$, то для любых $N \in M,~v\in N$ справедлива формула\\
$D_k(v)\equiv 0\mod{{\bf Z}[F]\cdot(N-1)},~k\in J\setminus K\Rightarrow v\in [N,N].$\\
Докажем обратное утверждение. Предположим противное. Пусть $1\neq v\in H\cap R$.\\
Так как $D_k(v)= 0,~k\in J\setminus K$, то $D_k(v)\equiv 0\mod{{\bf Z}[F]\cdot(R-1)},~k\in J\setminus K$. Следовательно, $v\in [R,R]=R^{(1)}$.\\
Так как $D_k(v)= 0,~k\in J\setminus K$, то $D_k(v)\equiv 0\mod{{\bf Z}[F]\cdot(R^{(1)}-1)},~k\in J\setminus K$. Следовательно, $v\in [R^{(1)},R^{(1)}]=R^{(2)}$.
Продолжая аналогичные рассуждения, получим
$v\in\bigcap\limits_{n \in {\bf N}} R^{(n)}=1$ --- противоречие.
\end{proof}

\section{О про-$p$-группах с одним определяющим соотношением}

\noindent Пусть $F_n$ --- дискретная свободная группа конечного ранга $n$, $\{x_1,\ldots ,x_n\}$ --- ее свободные
порождающие, $N_j$, $j\in J$ --- семейство всех нормальных подгрупп группы $F_n$,
индекс которых равен степени простого числа $p$, $I_{jm}$ --- идеал кольца ${\bf Z}[F_n]$, порожденный элементами $(p^m,1-h\mid h\in N_j,m\in {\bf N})$, ${\bf Z}_p$ --- кольцо целых $p$-адических чисел.\\
Обозначим через\\
$F(n)$ --- свободную про-$p$-группу с системой образующих $\{x_1,\ldots ,x_n\}$,
являющуюся пополнением группы $F_n$ в топологии, определяемой семейством $N_j$, $j\in J$,\\
через $\{U_n\mid n\in {\bf N}\}$ --- счетную убывающую систему открытых нормальных подгрупп группы $F(n)$, составляющих базу окрестностей единицы,\\
через ${\bf Z}_p[[F(n)]]=\varprojlim {\bf Z}_p[F/N_j]$ --- пополненное групповое кольцо свободной про-$p$-группы $F(n)$,
являющееся пополнением кольца ${\bf Z}[F_n]$ в топологии, определяемой семейством $I_{j,m}$, $j\in J,m\in {\bf N}$.\\
Производные Фокса $D_k$ кольца ${\bf Z}[F_n]$ и отображение $\varepsilon:{\bf Z}[F_n]\to {\bf Z}$ непрерывны в топологии
кольца ${\bf Z}[F_n]$, определяемой семейством идеалов $I_{j,m}$,
кольцо ${\bf Z}[F_n]$ --- всюду плотное подмножество в компактном кольце
${\bf Z}_p[[F(n)]]$, следовательно, существуют непрерывные продолжения на ${\bf Z}_p[[F(n)]]$ отображений $D_k$ и $\varepsilon$.\\ Обозначим эти продолжения через
$\partial_k$ и $\bar{\varepsilon}$ соответственно.
Из непрерывности отображений $\partial_k$ и $\bar{\varepsilon}$ следуют равенства:
\begin{gather}
\partial_k(v+w)=\partial_k(v)+\partial_k(w),\notag\\
\partial_k(vw)=\partial_k(v)w+\bar{\varepsilon}(v)\partial_k(w),\notag
\end{gather}
для любых $v,~w\in {\bf Z}_p[[F(n)]]$, $k = 1,\ldots ,n$.\\[3pt]
Доказательство теоремы \ref{tm1_gr}.
Рассмотрим $\{v_m\}$ --- последовательность элементов группы $F_n$, такую, что  $vv_m^{-1}\in U_m$.
Полагаем $\mathfrak{N}_{jt}$ --- идеал кольца ${\bf Z}_p[[F(n)]]$, порожденный элементами $(p^t,1-h\mid h\in NU_j,t\in {\bf N})$.
Ввиду формулы (\ref{end_gr}), для любого $\mathfrak{N}_{jt}$ найдется номер $L_{jt}$
такой, что
\begin{eqnarray}\label{1tm02_pr_gr}
D_k(v_l)\equiv ~0\mod{\mathfrak{N}_{jt}},~k\in J\setminus K,~l>L_{jt}.
\end{eqnarray}
Из (\ref{1tm02_pr_gr}), ввиду теоремы \ref{tm3_gr}, вытекает, что при $l>L_{jt}$ справедлива формула
\begin{eqnarray*}
v_l\equiv \hat{v}_l\prod\limits_{i=1}^{k_l}v_{l_i}^{f_{l_i}}\mod{[N,N]N^{p^t}U_j},
\end{eqnarray*}
где $\hat{v}_l\in F(K)$, $f_{l_i}\in F(n)$, $v_{l_i}\in F(K)\cap (NU_j)$.
Следовательно, для любых $t\in {\bf N}$ и $v_n$ найдутся $\hat{v}_n\in F(K)$, $\bar{v}_n\in (F(K)\cap N)^{F(n)}$ такие, что
\begin{eqnarray}\label{1tm02_pr_gr2}
v_n\equiv\hat{v}_n\bar{v}_n\mod{[N,N]N^{p^t}U_n}.
\end{eqnarray}
Ввиду (\ref{1tm02_pr_gr2}), $v$ принадлежит замыканиям в $F(n)$ групп
\begin{eqnarray*}
F(K)(F(K)\cap N)^{F(n)}[N,N]N^{p^t},~t=1,2,\ldots,
\end{eqnarray*}
следовательно, $v$ принадлежит замыканию $F(K)(F(K)\cap N)^{F(n)}[N,N]$ в $F(n)$.\\[3pt]
Доказательство следствия \ref{crl_pro_p_gr}. Предположим противное.\\
Тогда $H\cap R=1$ и $\pi(\partial_n(r))= 0$.
Теорема \ref{tm1_gr} показывает, что $r\in [R,R]$.\\
Следовательно, $r\in\bigcap\limits_{n \in {\bf N}} R^{(n)}=1$ --- противоречие.
\begin{corollary}\label{crl_pro_p_gr2}
Пусть $F(n)$ --- свободная про-$p$-группа с системой образующих $\{x_1,\ldots ,x_n\}$, $K\subset \{1,\ldots ,n\}$,
$H$ --- свободная про-$p$-группа с системой образующих $\{x_j \mid j\in K\}$, $R$ --- замкнутая нормальная подгруппа в $F(n)$, $\{\partial_1,\ldots ,\partial_n\}$ --- производные Фокса алгебры ${\bf Z}_p[[F(n)]]$, $M$ --- множество всех тех замкнутых нормальных подгрупп группы $F(n)$, которые принадлежат $R$.
Если (и только если) $H\cap R=1$, то для любых $N \in M,~v\in N$ справедлива формула
\begin{eqnarray*}
\pi_N(\partial_k(v))= 0,~k\in \{1,\ldots ,n\}\setminus K\Rightarrow v\in [N,N],
\end{eqnarray*}
где $\pi_N:{\bf Z}_p[[F(n)]]\rightarrow {\bf Z}_p[[F(n)/N]]$ --- гомоморфизм, индуцированный
естественным гомоморфизмом $F(n)\rightarrow F(n)/N$.
\end{corollary}
\begin{proof}
Теорема \ref{tm1_gr} показывает, что если $H\cap R=1$, то для любых $N \in M,~v\in N$ справедлива формула
\begin{eqnarray*}
\pi_N(\partial_k(v))= 0,~k\in \{1,\ldots ,n\}\setminus K\Rightarrow v\in [N,N].
\end{eqnarray*}
Докажем обратное утверждение. Предположим противное.\\
Пусть $1\neq v\in H\cap R$. Обозначим $\pi_j:{\bf Z}_p[[F(n)]]\rightarrow {\bf Z}_p[[F(n)/R^{(j)}]]$ --- гомоморфизм, индуцированный
естественным гомоморфизмом $F(n)\rightarrow F(n)/R^{(j)}$. Так как $\partial_k(v)= 0,~k\in \{1,\ldots ,n\}\setminus K$, то $\pi_1(\partial_k(v))= 0,~k\in \{1,\ldots ,n\}\setminus K$. Следовательно, $v\in [R,R]=R^{(1)}$.
Так как $\partial_k(v)= 0,~k\in \{1,\ldots ,n\}\setminus K$, то $\pi_2(\partial_k(v))= 0,~k\in \{1,\ldots ,n\}\setminus K$. Следовательно, $v\in [R^{(1)},R^{(1)}]=R^{(2)}$.
Продолжая аналогичные рассуждения, получим
$v\in\bigcap\limits_{n \in {\bf N}} R^{(n)}=1$ --- противоречие.
\end{proof}

\end{document}